\documentclass[graybox]{svmult}

\usepackage{type1cm}        % activate if the above 3 fonts are
\usepackage{makeidx}         % allows index generation
\usepackage{graphicx}        % standard LaTeX graphics tool
\usepackage{multicol}        % used for the two-column index
\usepackage[bottom]{footmisc}% places footnotes at page bottom

\usepackage{newtxtext}       % 
\usepackage{newtxmath}       % selects Times Roman as basic font

\makeindex             % used for the subject index
\newcommand{\R}{\mathbb{R}}
\newcommand{\C}{\mathbb{C}}
\newcommand{\N}{\mathbb N}

\DeclareMathOperator{\tr}{tr} 

\DeclareMathOperator{\diag}{diag}

\DeclareMathOperator{\spec}{spec}

\newcommand{\abs}[1]{\lvert#1\rvert}
\newcommand{\norm}[1]{\lVert#1\rVert}

\newcommand{\xii}{\abs{\xi}}

\begin{document}

\title*{Decay estimates for a Klein--Gordon model with time-periodic coefficients}
\author{Giovanni Girardi and Jens Wirth}
\institute{Giovanni Girardi \at Department of Mathematics, University of Bari, Via E. Orabona 4, 70125 Bari, ITALY, \email{giovanni.girardi@uniba.it}
\and Jens Wirth \at Department of Mathematics, University of Stuttgart, Pfaffenwaldring 57, 70569 Stuttgart, GERMANY,  \email{jens.wirth@mathematik.uni-stuttgart.de}}
\maketitle

\abstract*{In this paper we consider a Klein--Gordon model with time-dependent periodic coefficients. The aim is to investigate how the presence of the mass term influences energy estimates with respect to the case of vanishing mass, already treated in \cite{Wirth2008}. The approach is based on a diagonalisation argument for high frequencies and a contradiction argument for bounded frequencies.}

\abstract{In this paper we consider a Klein--Gordon model with time-dependent periodic coefficients. The aim is to investigate how the presence of the mass term influences energy estimates with respect to the case of vanishing mass, already treated in \cite{Wirth2008}. The approach is based on a diagonalisation argument for high frequencies and a contradiction argument for bounded frequencies.}

\section{Introduction}

In \cite{Wirth2008} the second author considered the linear Cauchy problem for a damped wave equation with time-periodic dissipation term 
$b(t)$,
\begin{equation}
\label{eq:CPWirth}
\begin{cases} 
u_{tt}-\Delta u + 2b(t)u_t=0,\\
u(0,x)=u_0(x), \quad u_t(0,x)=u_1(x),
\end{cases}
\end{equation}
and proved that the solution to \eqref{eq:CPWirth} satisfies the well-known Matsumura-type estimate obtained for constant dissipation by A.~Matsumura in \cite{Mat76}, that is 
\begin{equation}
\label{eq:matzumuraestimate}
\norm{\partial_t^k\nabla^j u(t,\cdot)}_{L^2}\leq C (1+t)^{-j-\frac{k}{2}}(\|u_0\|_{H^{j+k}}+\|u_1\|_{H^{j+k-1}}),
\end{equation}
for $j,\,k =0,\,1$ and $C$ a positive constant independent on the initial data. In this paper we generalise these results and consider the Cauchy problem 
\begin{equation}
\label{eq:CPmain}
\begin{cases} 
u_{tt}-\Delta u + 2b(t)u_t+m^2(t)u=0,\\
u(0,x)=u_0(x), \quad u_t(0,x)=u_1(x)
\end{cases}
\end{equation}
with positive time-periodic dissipation $b(t)$ and mass $m(t)$. We study how the presence of a periodic mass term influences the decay estimates for the solution to \eqref{eq:CPmain}.

Let us first explain, why such a problem is interesting and how it relates to known results from the literature.
There exist many papers in which decay estimates for the solution to wave models of the form
\eqref{eq:CPmain} are investigated under different assumptions on the coefficients $b(t)$ and $m(t)$. 
The survey articles \cite{Rei04} and \cite{Wir10} provide for an overview of results; moreover, we refer to the works of M.~Reissig and K.~Yagdjian~\cite{ReiYa00}, of F.~Hirosawa and M.~Reissig~\cite{HiRe06}, of M.~Reissig and J.~Smith~\cite{ReiSm05}, as well as the papers of the second author~\cite{Wirth2006}, \cite{Wirth2007}. In the latter two papers a classification of dissipation terms as \textit{non-effective} or \textit{effective} is introduced, which distinguishes the dissipation terms according to their strength and influence on the large-time behaviour of solutions. In all these results a control on the amount of oscillations present in the coefficients is essential. 

To understand this and the meaning of this classification we consider the Cauchy problem \eqref{eq:CPWirth} with the coefficient $b$ assumed to be a bounded, non-negative, sufficiently smooth function satisfying a condition of the form
\begin{equation}
\label{eq:oscillations}
| \partial_t^k b(t)|\leq C_k \frac{b(t)}{(1+t)^k}\qquad \text{ for } k=1,\,2.
\end{equation}
Then, we distinguish between two cases. First, if 
\begin{equation} 
\limsup_{t\to\infty} tb(t)<1.
\end{equation}
we say that $b$ is \textit{non-effective}, in the sense that the solution behaves in an asymptotic sense like a free wave multiplied by a decay factor, that is there exists a solution $v=v(t,x)$ to the wave equation $v_{tt} - \Delta v = 0$ such that
\[ \begin{pmatrix}
\nabla u(t,x) \\ u_t(t,x)
\end{pmatrix} \sim \frac{1}{\lambda(t)}\begin{pmatrix}
\nabla v(t,x) \\ v_t(t,x) 
\end{pmatrix}, \qquad t\to \infty, \]
the asymptotic equivalence understood in an appropriate $L^p$-sense and with $\lambda=\lambda(t)$ given as
\[\lambda(t)=\exp\Big(\frac{1}{2}\int_0^t b(\tau)\,d\tau\Big).\]
The initial data to the free wave $v=v(t,x)$ are uniquely determined by the solution $u=u(t,x)$ and thus by the initial data $u_0$ and $u_1$. Thus, a modified form of scattering is valid. On the other hand, if
\[ \lim_{t\to \infty} tb(t)=\infty\]
holds true we say that the dissipation $b$ is \textit{effective}; in this case solutions to damped wave equation are asymptotically related to solutions $w= w(t,x)$ of the parabolic heat equation $w_t=\Delta w$, i.e.
\[ u(t,x)\sim w(t,x)\]
holds true again in an appropriate $L^p$-sense. This can be made precise in the form of the so-called \textit{diffusion phenomenon} for damped waves; see \cite{Wirth07} for the time-dependent dissipation case or the papers of Nishihara \cite{Nish03} and Narazaki \cite{Nar04} for the case of constant dissipation. 

Wave models with mass and dissipation of the form \eqref{eq:CPmain} were considered by the second author and Nunes in \cite{NunesWirth}. This paper provides in particular $L^p-L^q$ decay estimates in the non-effective case.
In \cite{DAGR18} the first author considered with M.~D'Abbico and M.~Reissig the Cauchy problem \eqref{eq:CPmain} in the case in which the damping term is effective and dominates the mass term, i.e. $m(t)=o(b(t))$ as $t\to \infty$, again under control assumptions on the oscillations of the coefficients. In that paper it is shown that under a simple condition on the interaction between~$b(t)$ and~$m(t)$, one can prove that the solutions to \eqref{eq:CPmain} satisifies the estimate 
\begin{align}
\norm{u(t,\cdot)}_{L^2}
\label{eq:estimateDAGR}   & \leq C\,\gamma(t)\,\norm{(u_0,u_1)}_{H^1\times L^2},\\
\intertext{where we define}
\gamma(t)
\label{eq:gamma}  & = \exp \left(-\int_0^t \frac{m^2(\tau)}{b(\tau)}\,d\tau \right).
\end{align}
Thus, the decreasing function $\gamma=\gamma(t)$ in~\eqref{eq:gamma} represents the influence on the estimates of the mass term with respect to the damping term. In particular, estimate \eqref{eq:estimateDAGR} shows that the presence of the mass term produces an additional decay which becomes faster as the mass term becomes more influent. In fact, in \cite{GIRNTNP2019} the first author proved an exponential decay in the case of dominant mass, that is 
\begin{align}
\norm{u(t,\cdot)}_{L^2}
\label{eq:estimateGir2019}   & \leq C\,\exp \left(-\delta\int_0^t b(\tau)\,d\tau \right)\,\norm{(u_0,u_1)}_{H^1\times L^2},
\end{align}
provided that $\liminf_{t\to\infty} m(t)/b(t) > 1/4$. 
This latter estimate is almost the same as for the solution to the Cauchy problem to a damped Klein--Gordon model with constant coefficients $b(t)\equiv 1$ and $m(t)\equiv 1$, that is
\begin{equation}
\label{dampedKG}
\begin{cases}
u_{tt}-\Delta u+u_t+u=0,\\
u(0,x)=u_0(x), \quad u_t(0,x)=u_1(x).
\end{cases}
\end{equation}

All these cited papers have in common that they use assumptions on derivatives of the coefficients as in \eqref{eq:oscillations} 
to avoid a bad influence of oscillations. That oscillations may have deteriorating influences was shown
for example by K.~Yagdjian in \cite{Yag01} for a wave equation with time-periodic speed of
propagation. In this case (many) solutions have exponentially growing energy. Controlling oscillations is done by requiring estimates for derivatives of the coefficients. 

It is clear that for dissipative wave equations oscillations in the positive dissipation term can not lead to solutions with increasing energy. Therefore, it is interesting to ask whether conditions on derivatives of the coefficient are indeed necessary for proving large-time decay estimates for solutions of \eqref{eq:CPWirth}. A first step to look into that was done in \cite{Wirth2008}, where the author proved that the solution to \eqref{eq:CPWirth} satisifies estimate \eqref{eq:matzumuraestimate} without any condition on the oscillations of $b=b(t)$ provided that $b$ is periodic. This led to the conjecture that estimate \eqref{eq:matzumuraestimate} can be obtained with a general dissipation term $b=b(t)$, with $tb(t)\to\infty$, without further assumptions on derivatives. However, it is still an open problem how to prove such a result.

In the present paper we also avoid assumptions on the derivatives of the coefficients $b(t)$ and $m(t)$ assuming only that they are positive, periodic and of bounded variation. We are going to prove an exponential decay by using the same technique used as in \cite{Wirth2008} combined with a perturbation argument for the mass term. We remark that the presence of the mass term simplifies the study of the estimates at small frequencies; in fact, in this zone it is not necessary to use tools of Floquet theory as in the case of vanishing mass: we use only a contradiction argument together with some results of spectral theory of matrices.

The study of decay estimates for the solution to the linear problem \eqref{eq:CPmain} has an important application in the study of global (in time) existence results for the corresponding nonlinear problem  
\begin{equation*}
\begin{cases} 
u_{tt}-\Delta u + 2b(t)u_t+m^2(t)u=h(t,u),\\
u(0,x)=u_0(x), \quad u_t(0,x)=u_1(x)
\end{cases}
\end{equation*}
with nonlinearity $h(t,u)=(1+\int_0^t 1/b(\tau)\,d\tau)^\omega |u(t,\cdot)|^p$ for a $\omega \in[-1,\infty)$. Such applications can be found 
for example in \cite{DA13}, \cite{DALR13} in the purely dissipative case and in \cite{DAGR18,GIR2019,GIRNTNP2019} for equations including mass terms. 

The paper is organized as follows: In Section \ref{sec:main results} we give the basic assumptions on the Cauchy problem and we state our main results that are Theorem \ref{th:Thconstant} and Theorem \ref{th:Thperturbed}; in Section \ref{sec:representationsolution} we make considerations and discuss properties of the fundamental solution to \eqref{eq:CPmain} and the associated monodromy operator. In Section \ref{sec:smallfreqconstant} we treat the case of constant mass for small frequencies and we prove a fundamental lemma useful for the proof of the main theorems. Finally in Section \ref{sec:Proofofresults} the main theorems are proved.

\section{Main results}\label{sec:main results}
In this paper we suppose that the coefficient $b=b(t)$ is a non-negative and continuous periodic function of bounded variation, i.e., we assume that its weak derivative is essentially bounded, $b'\in L^{\infty}$. We further suppose that the coefficient $m=m(t)$ is measurable and periodic with the same period. We denote the period of both coefficients by $T$. The first result concerns constant mass terms and provides an exponential decay result.

\begin{theorem}
	\label{th:Thconstant}
	Suppose $m\equiv m_0\in\R$ is constant. There exists $\delta>0$ such that the solution $u=u(t,x)$ to the Cauchy problem \eqref{eq:CPmain} satisfies 
	\begin{align*}
	\| u(t,\cdot)\|_{L^2}&\leq Ce^{-\delta t}(\|u_0\|_{L^2}+\|u_1\|_{H^{-1}}),\\
	\| \nabla u(t,\cdot)\|_{L^2}&\leq Ce^{-\delta t}(\|u_0\|_{H^1}+\|u_1\|_{L^2}),\\ 
	\|u_t(t,\cdot)\|_{L^2}&\leq Ce^{-\delta t}(\|u_0\|_{H^1}+\|u_1\|_{L^2}),
	\end{align*}
	where $\delta$ and $C$ are positive constants depending on the coefficient $b$ and on $m_0$.
\end{theorem}

If the mass term is  non-constant, the exponential decay is obtained under a smallness condition for the deviation of the mass-term from a constant. 

\begin{theorem}
	\label{th:Thperturbed}
	Let $m_0\in \R$ and $m_1=m_1(t)$ a measurable $T$-periodic function such that $\sup_{t\geq 0} |m_1(t)|=1$. Then, there exists $\epsilon$ sufficiently small such that the solution to 
	\begin{equation}
	\label{eq:CPepsilon}
	\begin{cases} 
	u_{tt}-\Delta u + 2b(t)u_t+m_\epsilon^2(t)u=0,\\
	u(0,x)=u_0(x), \quad u_t(0,x)=u_1(x)
	\end{cases}
	\end{equation}
	with $m_\epsilon^2(t)=m_0^2+\epsilon m_1(t)$  satisfies
	\begin{align*}
	\| u(t,\cdot)\|_{L^2}&\leq Ce^{-\sigma t}(\|u_0\|_{L^2}+\|u_1\|_{H^{-1}}),\\
	\| \nabla u(t,\cdot)\|_{L^2}&\leq Ce^{-\sigma t}(\|u_0\|_{H^1}+\|u_1\|_{L^2}),\\ 
	\|u_t(t,\cdot)\|_{L^2}&\leq Ce^{-\sigma t}(\|u_0\|_{H^1}+\|u_1\|_{L^2}),
	\end{align*}
	where $\sigma$ and $C$ are positive constant depending on $m_0,\, m_1$, $b$ and $\epsilon$.
\end{theorem}
\begin{remark}
It is still an open problem to understand which is the largest value that $\epsilon$ can assume in order to guarantee an exponential decay of the energy. A possible estimate of $\epsilon $ is given in the proof of Theorem \ref{th:Thperturbed}: from estimate \eqref{eq:epsilon} it is clear that the value of $\epsilon$ depends on how large we choose $N$, such that the line $|\xi|=N$ divides the phase space in small and large frequencies. In particular, the value of $N$ depends only on the dissipation and does not depend on the mass term.
\end{remark}

\section{Representation of solution} \label{sec:representationsolution}
In a  first step we derive properties of the representation of solutions for the Cauchy problem
\begin{equation}
\label{eq:CPconstant}
u_{tt}-\Delta u + 2b(t)u_t+m^2(t)u=0, \qquad u(0,x)=u_0(x), \qquad u_t(0,x)=u_1(x),
\end{equation}
with $b=b(t)\geq 0$ and $m=m(t)\geq 0$ both periodic of period $T$. We denote the mean value of $b(t)$ as
\[ \beta=\frac{1}{T}\int_0^T b(t)\,dt. \]
A partial Fourier transform with respect to the spatial variables reduces the problem to an ordinary differential equation 
\begin{equation}
\label{eq:fouriertransformed}
\hat{u}_{tt}+\xii^2 \hat{u}+2b(t)\hat{u}_t+m^2(t)\hat{u}=0, 
\end{equation}
parameterised by $\xii\in\R$. To reformulate this as first order system, we introduce the symbol $\langle \xi \rangle_{m(t)}:=\sqrt{\xii^2+m^2(t)}\,$ and we define the new variable $V=(\langle \xi \rangle_{m(t)} \hat{u}, D_t \hat{u})^T$. Then
we obtain the system $D_tV=A(t,\xi)V$ with
\begin{equation}
\label{eq:mainsystem}
A(t,\xi)= \begin{pmatrix}
0 &\langle \xi \rangle_{m(t)} \\ \langle \xi \rangle_{m(t)}& 2ib(t)
\end{pmatrix},
\end{equation}
using the Fourier derivative $D_t=-i\partial_t$. We want to study the fundamental solution $\mathcal{E}=\mathcal{E}(t,s,\xi)$ to \eqref{eq:mainsystem}, that is the matrix-valued solution to the Cauchy problem
\begin{equation} 
\label{eq:system}
D_t \mathcal{E}(t,s,\xi)=A(t,\xi)\mathcal{E}(t,s,\xi), \qquad \mathcal{E}(s,s,\xi)=I. 
\end{equation}
In particular, we consider the family of monodromy matrices $\mathcal{M}(t,\xi)=\mathcal{E}(t+T,t,\xi)$. 
The fundamental solution to \eqref{eq:system} can be represented by the Peano--Baker series
\begin{equation}
\label{eq:Erepresntation}
\mathcal{E}(t,s,\xi)=I+\sum_{\ell=1}^\infty  i^k \int_s^tA(t_1,\xi)\int_s^{t_1}A(t_2,\xi)\cdots \int_s^{t_{\ell-1}}A(t_\ell,\xi)\,dt_\ell\cdots dt_1.
\end{equation}
The $T$-periodicity of coefficients implies periodicity of the matrix $A(t,\xi)$ and hence the $T$-translation invariance of the fundamental solution, i.e. $\mathcal{E}(t+T,s+T,\xi)=\mathcal{E}(t,s,\xi)$. Thus, the the monodromy matrix $\mathcal{M}(t,\xi)$ is $T$-periodic. 
Moreover, since $\mathcal E(t,s,\xi)\mathcal E(s,t,\xi) = I$ it follows that $\mathcal{E}(t,s,\xi)$ satisfies $D_s \mathcal{E}(t,s,\xi)=-\mathcal{E}(t,s,\xi)A(s,\xi)$, and, therefore, $\mathcal{M}(t,\xi)$ satisfies the equation
\begin{equation*}
D_t\mathcal{M}(t,\xi)=[A(t,\xi), \mathcal{M}(t,\xi)], \qquad  \mathcal{M}(T,\xi)=  \mathcal{M}(0,\xi).
\end{equation*}
In what follows we will distinguish between small and large frequencies and provide estimates for $\mathcal M$.
\subsection{Large Frequencies} \label{sec:largefrequencies}
For large frequencies we want to prove that the monodromy matrix is uniformly contractive, i.e.
\begin{equation}
\label{eq:contractive}
\|\mathcal{M}(t,\xi)\|<1
\end{equation}
holds true uniformly in $t\in [0,T]$ and $|\xi|\ge N$ for a constant $N$ chosen large enough. The choice of $N$ does not depend on the coefficient $m=m(t)$. In order to prove \eqref{eq:contractive} we apply two steps of diagonalization. We consider the unitary matrices 
\[ M= \frac{1}{\sqrt{2}}\begin{pmatrix}
1&-1\\1& 1
\end{pmatrix} \qquad M^{-1}=\frac{1}{\sqrt{2}}\begin{pmatrix}
1&1\\-1& 1
\end{pmatrix} \]
and define the new variable $V^{(0)}=M^{-1}V$, which satisfies
\[D_tV^{(0)}= (D(t,\xi)+R(t,\xi))V^{(0)} \] 
with 
\[ D(t,	\xi)=\begin{pmatrix}
\langle \xi \rangle_{m(t)} & 0 \\
0 & -\langle \xi \rangle_{m(t)}
\end{pmatrix},\qquad 
R(t,\xi)=ib(t)\begin{pmatrix}
1 & 1 \\
1 & 1
\end{pmatrix}.\]
Next, we define $D_1=D+\diag R$ and $R_1=R-\diag R$ and construct a matrix $N_1=N_1(t,\xi)$ with
\begin{equation}
\label{eq:N1}
D_tN_1= [D_1,N_1]+R_1,
\end{equation}
and $N_1(0,\xi)=I$.
Thus, the requirement for $N_1$ is equivalent to the operator identity 
\[ (D_t-D_1-R_1)N_1-N_1(D_t-D_1)= D_t N_1- [D_1,N_1]-R_1N_1=R_1(I-N_1).\]
Hence by denoting  $R_2=-N_1^{-1}R_1(I-N_1)$ we obtain
\[ (D_t-D_1-R_1)N_1= N_1(D_t-D_1-R_2) \] 
and as a consequence, provided that $N_1$ is invertible, we obtain that the new unknown $V^{(1)}=N_1^{-1}V^{(0)}$ satisfies the transformed equation
\[ D_tV^{(1)}= (D_1+R_2)V^{(1)}\]
with improved remainder allowing us later on to prove \eqref{eq:contractive}. 

Since $N_1=N_1(t,\xi)$ satisfies equation \eqref{eq:N1} and $D_1$ is diagonal, we find $D_t\diag N_1=0$. Thus, we can use a matrix $N_1$ of the form
\[N_1=\begin{pmatrix}
1 & n^-\\ n^+ & 1
\end{pmatrix},\]
with
\[ D_t n^\pm(t,\xi)= \mp\langle \xi \rangle_{m(t)}n^\pm(t,\xi)+i b(t).\]
The initial conditions $n^\pm(0,\xi)=0$ giving $N_1(0,\xi)=I$ imply
\[ n^\pm(t,\xi)=\int_0^t e^{\mp i\int_s^t\langle \xi \rangle_{m(r)}\,dr}b(s)\,ds. \]
Integrating by parts, we obtain
\[ |n^\pm(t,\xi)|=\Big|\Big[ \frac{\mp i}{\langle \xi \rangle_{m(s)}}e^{\mp i\int_s^t\langle \xi \rangle_{m(r)}\,dr}b(s)\Big]_0^t-\int_0^t \frac{\mp i}{\langle \xi \rangle_{m(s)}}e^{\mp i\int_s^t\langle \xi \rangle_{m(r)}\,dr}b'(s)ds \Big| \]
and using that $b=b(t)$ is of bounded variation we find a constant $C>0$ such that 
\[ |n^\pm(t,\xi)|\leq C(1+t)|\xi|^{-1}.\]
Thus we get that $n^\pm(t,\xi)\to 0$ when $\xi\to \infty$, uniformly in $[0,2T]$. Then we can conclude that $N_1(t,\xi)\to I$ and therefore $N^{-1}(t,\xi)\to I$ uniformly in $t\in [0,2T]$ as $|\xi|\to\infty$. Hence $\|R_2(t,\xi)\|\to 0$ as $|\xi|\to\infty$ uniformly in $t\in[0,2T]$.  Thus the supremum on the left hand side in the following formula tends to $1$ as $N\to\infty$ and we fix $N$ such that
\begin{equation}
\label{eq:suplargefreq}
\sup_{|\xi|\ge N} \sup_{t\in [0,T]} \| N_1(t+T,\xi)\| e^{\int_t^{t+T}\|R_2(s,\xi)\|ds}\|N_1^{-1}(t,\xi)\|\leq e^{\beta T/2}
\end{equation}
holds true. Note that this choice of $N$ can be made indepent of the coefficient $m=m(t)$, in fact $R_1=R_1(t,\xi)$ does not depend on $m(t)$, and $N_1(t,\xi),\, N_1^{-1}(t,\xi)$ tend both to $I$ uniformly with respect to $m(t)$. 

In order to prove the desired estimate $\eqref{eq:contractive}$ we go back to the original problem. We define $\lambda(t):=\exp( \int_0^t b(\tau)d\tau)$. Then, for each $|\xi|>N$ the fundamental solution $\mathcal{E}(t,s,\xi)$ to $D_tV=A(t,\xi)V$ with $A$ defined in \eqref{eq:mainsystem} is given by 
\begin{equation}
\label{eq:Erepresentation}
\mathcal{E}(t,s,\xi)=\frac{\lambda(s)}{\lambda(t)}M N_1(t,\xi)\tilde{\mathcal{E}}_0(t,s,\xi)Q(t,s,\xi)N_1^{-1}(t,\xi)M^{-1},
\end{equation}
for all $t\in [0,T]$, where 
\[ \tilde{\mathcal{E}}_0(t,s,\xi)= \begin{pmatrix}
e^{i\int_s^t \langle \xi \rangle_{m(\tau)}d\tau} & 0 \\0 & e^{-i\int_s^t \langle \xi \rangle_{m(\tau)}d\tau}
\end{pmatrix}\]
and $Q=Q(t,s,\xi)$ is the solution to the Cauchy problem
\[D_tQ(t,s,\xi)=\tilde{\mathcal{E}}_0(s,t,\xi)R_2(t,\xi)\tilde{\mathcal{E}}_0(t,s,\xi)Q(t,s,\xi), \qquad Q(s,s,\xi)=I.\]
Let $ \mathcal{R}_2(t,s,\xi)=\tilde{\mathcal{E}}_0(s,t,\xi)R_2(t,\xi)\tilde{\mathcal{E}}_0(t,s,\xi)$. Then by using the Peano-Beaker formula again we can represent $Q(t,s,\xi)$ as 
\[ Q(t,s,\xi)=I+\sum_{\ell=1}^\infty i^k \int_s^t \mathcal{R}_2(t_1,s,\xi)\int_s^{t_1} \mathcal{R}_2(t_2,t_1,\xi)\cdots \int_s^{t_{k-1}}  \mathcal{R}_2(t_k,t_{k-1},\xi) dt_k\dots dt_1. \] 
Since $\| \mathcal{R}_2(t,s,\xi)\|=\|R_2(t,\xi)\|$ we conclude 
\begin{equation}
\label{eq:Qestimate}
\|Q(t,s,\xi)\|\leq \exp\Big( \int_s^t \|R_2(\tau,\xi)\| d\tau \Big).
\end{equation} 
By \eqref{eq:Erepresentation} we can represent the monodromy matrix $\mathcal{M}(t,\xi)=\mathcal{E}(t,s,\xi)$ as
\[ \mathcal{M}(t,\xi)= \frac{\lambda(t)}{\lambda(t+T)}M N_1(t+T,\xi)\tilde{\mathcal{E}}_0(t+T,t,\xi)Q(t+T,t,\xi)N_1^{-1}(t+T,\xi)M^{-1}.\]
Since $\lambda(t)/\lambda(t+T)= e^{-\beta T}$ the desired result $\|\mathcal{M}(t,\xi)\|\le e^{-\beta T/2}<1$ for each $t\in [0,T]$ and each $|\xi|\geq N$ follows by \eqref{eq:suplargefreq} and \eqref{eq:Qestimate}. Hence we obtain

\begin{lemma}\label{lem1}
	There exists a constant $N$ depending only on $T$, $\|b'\|_\infty$ and $\|b\|_\infty$ such that the monodromy matrix $\mathcal M(t,\xi)$ satisfies
	\[
	\| \mathcal M(t,\xi) \| \le  e^{-\beta T/2}
	\]
	uniformly on $t\in\R$ and $|\xi|\ge N$ and independent of the mass term $m(t)$. 
\end{lemma}

\section{Small frequencies: constant mass} \label{sec:smallfreqconstant}
In this section we want to prove that there exists $k\in \N$ such that 
\begin{equation}
\label{eq:middlefreq}
\| \mathcal{M}^k(t,\xi)\|<1
\end{equation}  
uniformly in $ |\xi|\leq N$ and $t\in[0,T]$ provided that the mass term is constant.  Thus, in this section, we restrict our study to the Cauchy problem 
\begin{equation}
\label{eq:CPconstantmass}
v_{tt}-\Delta v +2b(t)v_t + m_0^2v=0
\qquad v(0,x)=v_0(x), \quad v_t(0,x)=v_1(x).
\end{equation}
In particular, we denote by $\mathcal{E}_0(t,s,\xi)$ the fundamental solution associated to the system $D_tV=A_0(t,\xi)V$ with
\begin{equation}
\label{eq:systemconstantmass}
A_0(t,\xi)= \begin{pmatrix}
0 &\langle \xi \rangle_{m_0} \\ \langle \xi \rangle_{m_0}& 2ib(t)
\end{pmatrix}, \qquad \langle \xi \rangle_{m_0}=\sqrt{|\xi|^2+m_0^2}.
\end{equation}
Let $\mathcal{M}_0(t,\xi)=\mathcal{E}_0(t+T,0,\xi)$ be the corresponding family of monodromy matrices. 
In order to get our aim we will prove at first that the spectrum $\spec \mathcal{M}_0(t,\xi)$ is contained in the open ball $\{ \eta\in \C \vert |\eta|<1 \}$.

Since it holds $$\mathcal{M}_0(t,\xi)\mathcal{E}_0(t,0,\xi)=\mathcal{E}_0(t+T,0,\xi)=\mathcal{E}_0(t+T,T,\xi)\mathcal{E}_0(T,0,\xi)=\mathcal{E}_0(t,0,\xi)\mathcal{M}_0(0,\xi),$$ we conclude that for each $t\in [0,T]$ the monodromy matrix $\mathcal{M}_0(t,\xi)$ is similar to $\mathcal{M}_0(0,\xi)$ and, hence,  has the same spectrum. Moreover, as both $b(t)$ and $m(t)$ are real; the  equation \eqref{eq:fouriertransformed} has real solutions and it follows that $\mathcal{M}_0(t,\xi)$ is similar to a real-valued matrix. Furthermore, by Liouville Theorem we know that
\begin{equation}
\det \mathcal{M}_0(0,\xi)= e^{i\int_0^T \tr A_0(\tau,\xi)d\tau }= e^{-2\beta T}.
\end{equation}
Hence, for each $\xi\in \R^n$ the eigenvalues $\eta_1(\xi),\,\eta_2(\xi)$ of $\mathcal{M}_0(0,\xi)$ are either real, in the form $\eta_2(\xi)= \eta_1^{-1}(\xi)e^{-2\beta T}$, or complex conjugate with $|\eta_1(\xi)|=|\eta_2(\xi)|=e^{-\beta T}$. In the latter case it is clear that $\spec \mathcal{M}_0(0,\xi)\subset \{ \xi\in \R^n | |\xi|=\exp(-\beta T)\}. $ In the case in which the eigenvalues are real we need to prove that for each $\xi \in \R^n$ both $\eta_1(\xi)$ and $\eta_2(\xi)$ have modulus less then $1$. We will prove this by using a contradiction argument.

Suppose that there exists $\bar{\xi}\in \R^n$ such that the monodromy matrix $\mathcal{M}_0(0,\bar{\xi})$ has an eigenvalue of modulus $1$, i.e, $\eta_1(\bar{\xi})= \pm 1$ and so $\eta_2(\bar{\xi})= \pm e^{-2\beta T}$. Let $\vec{c}=(c_1,c_2)$ be an eigenvector corresponding to $\eta_1(\bar{\xi})$. Then, we can find a domain $\Omega_R=\{ x\in \R^n | |x|\leq R \}$ (with $R$ depending on $m_0$) and a function $\Phi=\Phi(x)$ defined on $\Omega_R$ such that $-|\bar{\xi}|^2-m_0^2$ is an eigenvalue for the Dirichlet Laplacian with normal eigenfunction $\Phi=\Phi(x)$, i.e. 
\begin{equation}
\label{eq:laplacedirichlet}
-\Delta \Phi(x)= (|\bar{\xi}|^2+m_0^2)\Phi(x), \qquad \Phi(x)=0 \text{ on } \partial \Omega_R.
\end{equation}

Let us consider $v=v(t,x)$ the solution to the Cauchy problem, with Dirichlet boundary condition on $\Omega_R$ 
\begin{equation}
\label{eq:dampedwavedirichlet}
\begin{cases} v_{tt}-\Delta v+2b(t)v=0, \\  v(0,x)=c_1\langle \bar{\xi}\rangle_{m_0}^{-1}\Phi(x), \quad v_t(0,x)=ic_2 \Phi(x), \\ v(t,\cdot)\equiv 0 \quad \text{ on } \partial\Omega_R \text{ for each } t\geq 0. \end{cases}
\end{equation}
In particular, we look for a solution in the form
$$v(t,x)=f(t)\Phi(x),$$
and we show that $f=f(t)$ is $T$-periodic (or $2T$-periodic). Since, $\Phi=\Phi(x)$ satisfies the Dirichlet problem \eqref{eq:laplacedirichlet}, the partial differential equation $v_{tt}-\Delta v+2b(t)v=0$ turns into the ordinary differential equation $v_{tt}+|\bar{\xi}|^2 v+2b(t)v_t+m_0^2 v=0$, with $x$ regarded as a parameter. In particular, $f=f(t)$ satisfies the ordinary differential equation
\begin{equation}
\label{eq:odef}
f''(t)+2b(t)f'(t)+(|\bar\xi |^2+m_0^2)f(t)=0.
\end{equation}

Moreover, the corresponding solution $v(t,x)=f(t)\Phi(x),$ satisfies the Cauchy problem
\begin{align*}
D_t \begin{pmatrix}
\langle \bar\xi \rangle_{m_0} v(t,x)\\
D_t v(t,x)\end{pmatrix}&= \begin{pmatrix}
0 &\langle \bar\xi \rangle_{m_0} \\ \langle \bar\xi \rangle_{m_0}& 2ib(t)
\end{pmatrix}\begin{pmatrix}
\langle \bar\xi \rangle_{m_0} v(t,x)\\
D_t v(t,x)\end{pmatrix}\\ 
\begin{pmatrix}
\langle \bar\xi \rangle_{m_0} v(t,x)\\
D_t v(t,x)\end{pmatrix}\Big| _{t=0}&=  \begin{pmatrix}
c_1\\
c_2\end{pmatrix}\Phi(x).
\end{align*}
This system can be solved by using the fundamental solution $\mathcal{E}_0(t,0,\bar\xi)$; in particular, we have that 
\[ \begin{pmatrix}
\langle \bar\xi \rangle_{m_0} v(t,x)\\
D_t v(t,x)\end{pmatrix}\Big|_{t=T}= \mathcal{M}_0(0,\bar\xi)\begin{pmatrix}
c_1\\
c_2\end{pmatrix}\Phi(x)=\pm \begin{pmatrix}
c_1\\
c_2\end{pmatrix}\Phi(x). \]
We conclude that $f=f(t)$ is $T$-periodic (or $2T$-periodic) and $f(0)=c_1\langle \bar{\xi}\rangle_{m_0}^{-1}$. This gives a contradiction: if we denote the energy of this solution as 
\[
E(u,t)=1/2\|v_t(t,\cdot)\|_{L^2(\Omega_R)}^2+1/2\|\nabla v\|_{L^2(\Omega_R)}^2,\]
we obtain
\[ \frac{d}{dt}E(v,t)= -b(t)\|v_t\|_{L^2(\Omega_R)}^2=-b(t)|f'(t)|^2. \]
But, by integrating the previous equation we obtain that 
$$ -\int_0^T b(t)|f'(t)|^2 \,dt=0,$$
that is not possible since $f=f(t)$ can not be constant, by equation \eqref{eq:odef} as $(|\bar\xi|^2+m^2_0)>0$ for each $\bar\xi\in \R^n$. Thus, $\pm 1\notin \spec\mathcal{M}_0(t,\xi)$ for each $\xi\in \R^n$, and therefore the spectral radius $\rho(\mathcal{M}_0(t,\xi))<1$ for all $\xi \in \R^n$. \\
By the spectral radius formula, we know that $$\lim_{k\to \infty}\|\mathcal{M}_0^k(t,\xi)\|^{\frac{1}{k}}=\rho(\mathcal{M}_0(t,\xi))<1.$$ Thus, we conclude that for each $t\in [0,T]$ and $\xi\in\R^n$  there exists $k=k(t,\xi)\in \mathbb{N}$ such that 
\begin{equation}
\label{eq:Mk}
\|\mathcal{M}_0^k(t,\xi)\|<1.
\end{equation}
We want to show that we can find a number $k$ such that the condition \eqref{eq:Mk} holds uniformly with respect to $t\in [0,T]$ and $|\xi|\in [0,N]$.\\
Let us define for each $k\in \mathbb{N}$ the set $\mathcal{U}_k=\{(t,\xi)\in \R_+\times \R^n|\,\|\mathcal{M}_0^k(t,\xi)\|<1 \}$. It is open due to the continuity of the monodromy matrix ${M}_0^k(t,\xi)$; moreover, it holds $\mathcal{U}_k \subset \mathcal{U}_\ell$, for $k\leq \ell$. Then, by \eqref{eq:Mk} we have that the compact set $\mathcal{C}=\{(t,\xi)| 0\leq t\leq T, |\xi|\le N\}$ is contained in $\bigcup_{k}\mathcal{U}_k$. By compactness we find $k\in \mathbb{N}$ such that $\mathcal{C}\subset \mathcal{U}_k$. This concludes the proof of estimate \eqref{eq:middlefreq}. By continuity of $\mathcal M_0^k(t,\xi)$ in both variables, the estimate is uniform. Hence we obtain

\begin{lemma}\label{lem2}
	For constant mass term $m_0$ and fixed $N>0$ there exists a number $k$ such that the monodromy matrix for the problem with constant mass satisfies
	\[ \sup_{|\xi|\le N} \sup_{t\in [0,T]}\|\mathcal{M}_0^k(t,\xi)\|<1.\]
\end{lemma}

\section{Proof of the main theorems} \label{sec:Proofofresults}
\subsection{Proof of Theorem \ref{th:Thconstant}}
In order to prove Theorem \ref{th:Thconstant} we distinguish between small and large frequencies. 

Let $|\xi|\geq N$. Then the monodromy matrix $\mathcal{M}(t,\xi)$ is estimated in Lemma~\ref{lem1}. Let $t\geq 0$, $t=\ell T+s$, with $\ell\in \mathbb{N}$ and $s\in [0,T]$. Then, we obtain
\[ \|\mathcal{E}(t,0,\xi)\|=\|\mathcal{M}^\ell(s,\xi)\mathcal{E}(s,0,\xi)\|\leq  e^{-\ell \beta T/2} \|\mathcal{E}(s,0,\xi)\|. \]
Moreover, since $b(t)>0$ we know that $\|\mathcal{E}(s,0,\xi)\|\leq 1$ and therefore we find 
\[ \|\mathcal{E}(t,0,\xi)\|\leq e^{-\delta_0(t-T)},\]
by defining $\delta_0:=\beta/2>0$. We remark that this estimate for large frequencies is valid for arbitrary periodic mass terms.  

For the remainder of the proof assume that $m^2(t)\equiv m_0^2$ constant and $|\xi|\le N$. By Lemma~\ref{lem2} there exists $k\in \mathbb{N}$ depending only on $m_0$ such that the matrix $\mathcal{M}_0^k(t,\xi)$ is a contraction uniform in $t$ and $\xi$. Let $t=\ell k T+s \ge 0$ for some $\ell \in \N$ and $s\in [0,kT]$. Then, we obtain the exponential decay 
\begin{equation}
\label{eq:mainsmallconstant}
\|\mathcal{E}_0(t,0,\xi)\| = \| \mathcal M_0^{k\ell} (s,\xi) \mathcal E_0(s,0,\xi)\| \leq e^{-\delta_1(t-kT)},
\end{equation}
where we set $\delta_1:=(kT)^{-1}\log(c_1(N)^{-1})>0$ and 
\begin{equation}\label{c1def}
c_1(N):=\sup_{|\xi|\le N} \sup_{t\in [0,T]} \|\mathcal{M}_0^k(t,\xi)\|<1.
\end{equation}
Going back to the original problem \eqref{eq:CPmain}, we find
\[ \begin{pmatrix}
\langle \xi \rangle_{m_0} \hat{u}(t,\xi) \\D_t \hat{u}(t,\xi)
\end{pmatrix}= \mathcal{E}(t,0,\xi)\begin{pmatrix}
\langle \xi \rangle_{m_0} \hat{u}_0(\xi) \\ \hat{u}_1(\xi)
\end{pmatrix}.\]
Thus, we find that 
\begin{align*}
\| u(t,\cdot)\|_{L^2}&\leq \sup_{\xi\in\R^n}  \|\mathcal{E}_0(t,0,\xi)\|(\|u_0\|_{L^2}+\|u\|_{H^{-1}} ),\\
\| \nabla u(t,\cdot)\|_{L^2}&\leq \sup_{\xi\in\R^n} \|\mathcal{E}_0(t,0,\xi)\|(\|u_0\|_{H^1}+\|u_1\|_{L^2}),\\ 
\|u_t(t,\cdot)\|_{L^2}&\leq \sup_{\xi\in\R^n} \|\mathcal{E}_0(t,0,\xi)\|(\|u_0\|_{H^1}+\|u_1\|_{L^2}).
\end{align*}
The proof of Theorem \ref{th:Thconstant} with $C=e^{\delta_1 k T}$ follows immediately by estimate \eqref{eq:mainsmallconstant}.

\subsection{Proof of Theorem \ref{th:Thperturbed}}
Let $u=u(t,x)$ the solution to \eqref{eq:CPepsilon}
where $m_\epsilon^2(t)=m_0^2+\epsilon m_1(t)$, whit $ m_1(t)$ periodic of period $T$ and $m_0$ a sufficiently large constant such that $m_0^2+\epsilon m_1(t)>0$.
The corresponding system is 
\begin{equation}
\label{eq:perturbedsystem}
D_tV_\epsilon=A_\epsilon(t,\xi)V_\epsilon= \begin{pmatrix}
0 &\langle \xi \rangle_{m_\epsilon(t)} \\ \langle \xi \rangle_{m_\epsilon(t)}& 2ib(t)
\end{pmatrix}V_\epsilon,
\end{equation}
where $V_\epsilon=(\langle\xi\rangle_{ m_\epsilon(t)}\hat{u^\epsilon}, D_t \hat{u^\epsilon} )$.
In order to obtain our result we need to estimate $\| \mathcal{E}_\epsilon(t,0,\xi)\|$, where we denoted by $\mathcal{E}_\epsilon$ the fundamental solution to the system \eqref{eq:perturbedsystem}. In particular, $\mathcal{E}_0$ solves $D_tV_0=A_0(t,\xi)V_0$ where 
\begin{equation}
\label{eq:constantsystem}
D_tV_0=A_0(t,\xi)V_0= \begin{pmatrix}
0 &\langle \xi \rangle_{m_0} \\ \langle \xi \rangle_{m_0}& 2ib(t)
\end{pmatrix}V_0.
\end{equation}
We again distinguish between small and large frequencies. If $|\xi|\geq N$, as in the case of constant mass we conclude
$$\|\mathcal{E}_\epsilon(t,0,\xi)\|\leq e^{-\delta_0(t-T)},$$
where we recall $\delta_0=\beta/2>0$ by making use of Lemma~\ref{lem1}.

If $|\xi|\leq N$, there exists $k\in \mathbb{N}$ given by Lemma~\ref{lem2} such that the matrix $\mathcal{M}_0^k(t,\xi)$ is a contraction uniformly in $t\in [0,T]$ and $|\xi|\in [0,N]$. We write $t=\ell k T+s\ge0$ for some $\ell \in \mathbb N$ and $s\in [0,kT]$; then, we have
\begin{equation}
\label{eq:mainsmallperiodic}
\mathcal{E}_\epsilon(t,0,\xi) = \mathcal M_\epsilon^{k\ell} (s,\xi) \mathcal E_\epsilon(s,0,\xi);
\end{equation}
we can treat the fundamental solution as a perturbation of constant case
\[\begin{split}
\|\mathcal{E}_\epsilon(t,s,\xi)\|&\leq \|\mathcal{E}_\epsilon(t,s,\xi)-\mathcal{E}_0(t,s,\xi)\|+\|\mathcal{E}_0(t,s,\xi)\|\\&\leq \|\mathcal{E}_\epsilon(t,s,\xi)-\mathcal{E}_0(t,s,\xi)\|+ e^{-\delta(t-s-kT)},
\end{split}\]
where we recall $\delta_1=(kT)^{-1}\log(c_1(N)^{-1})>0$ and $c_1(N)$ as in \eqref{c1def}.
In order to estimate the difference  $\|\mathcal{E}_\epsilon(t,s,\xi)-\mathcal{E}_0(t,s,\xi)\|$ we use that for each $\epsilon\geq 0$ the fundamental solution $\mathcal{E}_\epsilon$ satisfies the integral equation
\[\mathcal{E}_\epsilon(t,s,\xi)=I+\int_s^t A_\epsilon(\tau,\xi) \mathcal{E}_\epsilon(\tau,s,\xi)\,ds,\]
such that
\begin{align*}
\mathcal{E}_\epsilon(t,s,\xi)-\mathcal{E}_0(t,s,\xi)=&\int_s^t A_\epsilon(\tau,\xi)(\mathcal{E}_\epsilon(\tau,s,\xi)-\mathcal{E}_0(\tau,s,\xi))\,ds\\&+\int_s^t (A_\epsilon(\tau,\xi)-A_0(\tau,\xi))\mathcal{E}_0(\tau,s,\xi)\,ds.
\end{align*}
By using the Gronwall inequality we get 
\[\|\mathcal{E}_\epsilon(t,s,\xi)-\mathcal{E}_0(t,s,\xi)\|\leq \int_s^t \|\mathcal{E}_0(\tau,s,\xi)\|\,\|A_\epsilon(\tau,\xi)-A_0(\tau,\xi)\|\,ds\cdot e^{\int_s^t \|A_\epsilon(\tau,\xi)\|\,d\tau}; \]
here, for any $\tau>0$ and $\xi\in \R^n$, since we are assuming $\displaystyle{\sup_{t\geq0} |m_1(t)|=1}$ we can estimate $$ \|A_\epsilon(\tau,\xi)-A_0(\tau,\xi)\|\leq \frac{\epsilon}{\langle\xi\rangle_{m_0}},\qquad \|A_0(\tau,\xi)\|\leq \langle \xi \rangle_{m_0}+2b(\tau),$$
and so $$ \|A_\epsilon(\tau,\xi)\|\leq C_\epsilon(\xi)+\langle \xi \rangle_{m_0}+2b(\tau), \qquad C_\epsilon(\xi)=\frac{\epsilon}{\langle\xi\rangle_{m_0}}.$$
Thus, recalling that $\mathcal{M}^k_\epsilon(s,\xi)= \mathcal{E}_\epsilon(s+kT,s,\xi)$, we find
\begin{align*} \| \mathcal{M}^k_\epsilon(s,\xi)-\mathcal{M}^k_0(s,\xi)\|&\leq C_\epsilon(\xi)e^{C_\epsilon(\xi) kT}e^{(\langle\xi\rangle_{m_0}+2\beta)k T}\int_s^{s+kT}\|\mathcal{E}_0(\tau,s,\xi)\|\,d\tau \\ &\leq C_\epsilon(\xi)e^{C_\epsilon(\xi) kT}e^{(\langle\xi\rangle_{m_0}+2\beta)k T}\int_s^{s+kT} e^{-\delta_1(\tau-s-kT)}\,d\tau \\
&\leq C_\epsilon(\xi)e^{C_\epsilon(\xi) kT}e^{(\langle\xi\rangle_{m_0}+2\beta)k T}\int_s^{s+kT} e^{-\delta_1(\tau-s-kT)}\,d\tau \\
&\leq \frac{C_\epsilon(\xi)}{\delta_1}e^{C_\epsilon(\xi) kT}e^{(\langle\xi\rangle_{m_0}+2\beta)k T}(e^{\delta_1 kT}-1).
\end{align*} 
Therefore, recalling that $\exp(\delta_1kT)=c_1(N)^{-1}$, we can conclude
\begin{align*}
\sup_{\xii\leq N}\sup_{s\in [0,T]}\| \mathcal{M}^k_\epsilon(s,\xi)\|&\leq \sup_{\xii\leq N}\sup_{s\in [0,T]}\|\mathcal{M}^k_0(s,\xi)\|\\
&\qquad + \sup_{\xii\leq N}\Big\{
\frac{C_\epsilon(\xi)}{\delta_1}e^{C_\epsilon(\xi) kT}e^{(\langle\xi\rangle_{m_0}+2\beta)k T}(c_1(N)^{-1}-1)\Big\}\\&=c_1(N)+ \sup_{\xii\leq N}\Big\{ \frac{C_\epsilon(\xi)}{\delta_1}e^{C_\epsilon(\xi) kT}e^{(\langle\xi\rangle_{m_0}+2\beta)k T}(c_1(N)^{-1}-1)\Big\}.
\end{align*}
By \eqref{eq:mainsmallperiodic} we get the desired result 
\[\sup_{\xii\leq N}\sup_{s\in [0,T]}\| \mathcal{M}^k_\epsilon(s,\xi)\|<1,\]
by choosing $\epsilon$ sufficiently small such that
\begin{equation}
\label{eq:epsilonrough}
 \frac{C_\epsilon(\xi)}{\delta_1}e^{C_\epsilon(\xi) kT}e^{(\langle\xi\rangle_{m_0}+2\beta)k T}(c_1(N)^{-1}-1)<1-c_1(N). 
\end{equation}
Let us introduce $W=W(x)$ the Lambert W-function defined in the set $\R_+:=\{x\in\R: x\geq 0\}$ such that for any $x\in \R_+$ it holds $x=W(x)e^{W(x)}$.  The function $W$ is increasing (see \cite{WLambert} for more details); thus, recalling the definition of $\delta_1$, we find that estimate \eqref{eq:epsilonrough} is equivalent to ask 
\[ \epsilon \leq \frac{\langle \xi \rangle_{m_0}}{kT}W\Big(c_1(N)\log(c_1(N)^{-1})e^{-(\langle\xi\rangle_{m_0}+\beta)kT}\Big),\]
for any $\xi \in [0,N]$, that is
\begin{equation}
\label{eq:epsilon}
 \epsilon \leq \frac{m_0}{kT}W\Big(c_1(N)\log(c_1(N)^{-1})e^{-(\langle N \rangle_{m_0}+\beta)kT}\Big).
\end{equation}

\begin{acknowledgement}
The paper is based on discussions the authors had during the stay of Giovanni Girardi at the University of Stuttgart in spring 2019. G.G. is grateful for the hospitality of the Department of Mathematics during his stay.
\end{acknowledgement}

%%%%%%%%%%%%%%%%%%%%%%%% referenc.tex %%%%%%%%%%%%%%%%%%%%%%%%%%%%%%
% sample references
% %
% Use this file as a template for your own input.
%
%%%%%%%%%%%%%%%%%%%%%%%% Springer-Verlag %%%%%%%%%%%%%%%%%%%%%%%%%%
%
% BibTeX users please use
% \bibliographystyle{}
% \bibliography{}
%

\end{document}